\documentclass[12pt]{amsart}

\usepackage{calc,pifont}
\usepackage{graphicx,adjustbox}
\usepackage{mathrsfs}
\usepackage{amssymb,amsbsy}
\usepackage{amscd}
\usepackage{enumerate}
\usepackage{setspace}
\usepackage{tikz}
\usepackage{bm}

\def\today{\number\day\space\ifcase\month\or
January\or February\or March\or April\or May\or June\or July\or
August\or September\or October\or November\or December\fi
\number\year}

\theoremstyle{definition}
\newtheorem{theorem}{Theorem}[section]
\newtheorem{lemma}[theorem]{Lemma}
\newtheorem{proposition}[theorem]{Proposition}
\newtheorem{definition}[theorem]{Definition}
\newtheorem{corollary}[theorem]{Corollary}

\newenvironment{pf}{{\emph{Proof:}}}{\QED}

\renewcommand{\phi}{\varphi}
\newcommand{\norma}[1]{\left\| #1 \right\|}
\newcommand{\inpr}[1]{\left\langle #1 \right\rangle}	
\newcommand{\conj}[1]{\left\{ #1 \right\}}
\newcommand{\ovl}[1]{\overline{#1}}	
\newcommand{\abs}[1]{\left| #1 \right|}
\newcommand{\parn}[1]{\left( #1 \right)}
\newcommand{\conv}[1]{\mathrm{conv}\parn{#1}}	
\newcommand{\orb}[1]{\mathrm{orb}\parn{#1}}
\newcommand{\conjz}{\conj{0,1}^\Z}	

\newcommand{\bz}{{\bf 0}}
\newcommand{\bo}{{\bf 1}}
\newcommand{\bes}{\begin{equation*}}
\newcommand{\ees}{\end{equation*}}
\newcommand{\beq}{\begin{equation}}
\newcommand{\eeq}{\end{equation}}
\newcommand{\im}{\texttt{i}}

\newcommand{\si}{\sigma}

\newcommand{\lm}{\lambda}

\newcommand{\Z}{{\mathbb{Z}}}
\newcommand{\R}{{\mathbb{R}}}
\newcommand{\C}{{\mathbb{C}}}

\newcommand{\D}{{\mathbb{D}}}

\newcommand{\QED}{\rule{0.4em}{2ex}}

\author{Christian Hern\'{a}ndez-Becerra}
\author{Benjam\'{\i}n A.~Itz\'{a}-Ortiz}

\address{Centro de Investigaci\'{o}n en Matem\'{a}ticas,
 Universidad Aut\'{o}noma del Estado de Hidalgo,
 Pachuca de Soto, Hidalgo, 42090, M\'{e}xico.}
\email[]{itza@uaeh.edu.mx}
\email[]{chriz.mate@gmail.com}

\thanks{}

\keywords{Tridiagonal operatores, pseudoergodic sequences, subshifts}

\begin{document}

\title[Tridiagonal operators for subshifts]{A class of tridiagonal operators associated to some subshifts}

\begin{abstract} 
 We consider a class of tridiagonal operators induced by not necessary pseudoergodic biinfinite sequences. Using only elementary techniques we prove that the numerical range of such operators is contained in the convex hull of the union of the numerical ranges of the operators corresponding to the  constant biinfinite sequences;  whilst the other inclusion is shown to hold when the constant sequences belong to the subshift generated by  the given biinfinite sequence. Applying recent results  by S.~N. Chandler-Wilde et~al. and R.~Hagger, which rely on limit operator techniques,  we are able  to provide more general results although the closure of the numerical range needs to be taken.
\end{abstract}

\maketitle

\section*{Introduction}

Let $b=(b_i)_{i\in\Z}$ be a biinfinite sequence in $\mathcal A^\Z$ where $\mathcal A$ is a finite set, called an alphabet. In this paper we study the operator $A_b\colon \ell^2(\Z)\to\ell^2(\Z)$  defined as the tridiagonal operator

 \[A_{b}= 
\left( \begin{array}{ccccccc}
\ddots & \ddots & & & & & \\
\ddots & 0 & 1 & & & & \\
& b_{-2} & 0 & 1 & & & \\
& & b_{-1} & \framebox[0.4cm][l]{0} & 1 & & \\
& & & b_{0} & 0 & 1 & \\
& & & & b_{1} & 0 & \ddots \\
& & & & & \ddots & \ddots
\end{array} \right) \]

\noindent where the rectangle marks the matrix entry at $(0,0)$. When the alphabet is the set $\{-1,1\}$,  the correspondig operators are related to the so called ``hopping sign model'' introduced in \cite{Feinberg} and subsequently studied in \cite{CWChLi10,CWChLi13,ChD,Hagger,HaggerJFA}.  We remark that  results found in the literature focus on the case when $b$ is a pseudoergodic sequence, that it to say,  when every finite sequence of $\pm 1$ appears somewhere in $b$. In this paper we aim to investigate the case when $b$ is not necessary pseudergodic and although some authors have in fact pointed out that some of their results hold if the pseudoergodic condition is dropped, we try to formalize this approach. For this reason, we begin by  thinking that $b$ is an element in a full shift space and as such it is more usual to consider the alphabet $\mathcal A$ to be the set $\{0,1\}$ rather that $\{-1,1\}$.  By means of elementary methods, we are able to establish  similar results to  the ones found in the literature and slighly generalize others, more concretely, we show that the numerical range of $A_b$, has an inclusion-wise maximal  but that such bound is sharp not only when $b$ is pseudoergodic but whenever the subshift generated by $b$ (the closure of the orbit of $b$) contains the constant biinfinite sequences. With this motivation we employ recent limit operator techniques to extend our results to more general alphabets; however, for these generalizations the closure of the numerical range must be taken. 

We feel our approach might be of interest as  it  employs elementary mathematics to obtain some results that are the motivation to state more gerenal results. We devide this work in three sections. In the first section we review some of the fundamental results needed in the rest of the paper, both operator theoretical concepts and symbolic dynamics notions. In the second section we apply elementary techniques to bound the numerical range of $A_b$, for any $b\in\{0,1\}\sp\Z$ and prove that such upper bound is sharp when the subshift generated by $b$ contains the constant sequences. In the last section we generalize  our results to more general alphabets by means of recent results which employ limit operator techniques. 

The authors gratefully acknowledge stimulating conversations with Rub\'en Martinez-Avenda\~no, Federico Menendez-Conde and Jorge Vi\-ve\-ros during the preparation of this paper.

\section{Background}
We review in this section some of the fundamental results and notation to be used throughout in the paper. 

\subsection{Operators}

  We are considering bounded linear operators on the Hilbert space $\ell^2(\Z)$. Hence, every time we refer to an operator we assume it is a bounded linear operator on $\ell^2(\Z)$. The inner product on $\ell^2(\Z)$ is denoted by $\inpr{\phantom{a},\phantom{b}}$.  For an operator $A$ we define the spectrum of $A$ as
	\bes
		\si(A)=\conj{\lm\in\C \ \middle | \ A-\lm I \mbox{ is not invertible}}
	\ees
The numerical range of $A$ is defined by
\[
   W(A)=\conj{ \inpr{Ax,x} \colon x\in \ell^2(\Z),\ \| x\|=1}
\]
 For easy reference, we state the following  well known   results (see e.g.\cite[Chapter~1]{RubMar07}) which are true in general Hilbert spaces.

\begin{theorem}\label{teospectrum}
	Let be $S$ and $T$ be operators and denote by $I$ the identity operator. Then
	\begin{enumerate}[(i)]
		\item $\si(T)\subset\C$, $\si(T)\neq\emptyset$ and $\si(T)$ is compact set.
			\item $W(I)=\conj{1}$, and for $\alpha,\beta\in\C$, $W(\alpha T + \beta I) = \alpha W(T) + \beta$.
			\item $W(T+S)\subseteq W(T) + W(S)$.
			\item $W(T)$ is a convex subset of $\C$.
			
              \item If $T$ is normal, then $\ovl{W(T)}$ (the closure of the numerical
range of $T$) is the convex hull of $\sigma(T)$.
             \item $\si(T)\subseteq\ovl{W(T)}$.
              \item If $T$ is self-adjoint then $W(T)\subset\R$ and $\si(T)\subset\R$.
	\end{enumerate}
\end{theorem}

\subsection{Symbolic dynamics}

For details about symbolic dynamics we refer the reader to \cite{LinMar95}.
An alphabet  $\mathcal A$ is nothing but a finite set.  We refer to  the elements of $\mathcal A$ as {symbols}. The set $\mathcal A\sp\Z$ of all biinfinite sequnces of symbols from $\mathcal A$ is termed the full shift on $\mathcal A$. 
When writing down an element in $\mathcal A\sp\Z$ it is costumary to distinguish the $0$th coordinate with a dot. For example a biinfinite sequence $b$  in $\conjz$ expressed as
\bes
b=(\ldots 0\,1\,1\,1\,0 \dot 1\,0\,0\,0\,1\,0\ldots)
\ees
so that its coordinates are $b_{-2}=1$, $b_{-1}=0$, $b_0=1$, $b_1=0$ and so on.

 Given a biinfinite sequence $b=(b_i)_{i\in\Z}$ in $\mathcal A\sp\Z$, a {block} (or {word}) of $b$ is defined to be a finite subsequence of $b$. The length of a block of $b$ is the number of symbols it contains. If $b$ belongs to $A^\Z$ and $i < j$, then we will denote the block of coordinates in $b$ from position $i$ to position $j$ by $b_{[i,j]}$. When a finite sequence $u$ of symbols of $\mathcal A$ satisfies that $u=b_{[i,j]}$ for some $i,j$, we will say that $u$ occurs in $b$. The blocks of the form $b_{[-k,k]}$ are called central blocks of $b$.
Note  then an equivalent way to say that $b\in\mathcal A$ is pseudoergodic is to require that  every block of every size of symbols in $\mathcal A$ occur in  $b$. 

The full shift $\mathcal A\sp\Z$ is actually a metric space with a metric given by
 	\bes
		\rho(b, c) = \left\{
		\begin{array}{cl}
			0				& \mbox{ if }\ b=c\\
			1				& \mbox{ if }\ b_0\neq c_0\\
			2^{-k}	& \mbox{ if $k$ is maximal so that $b_{[-k,k]} = c_{[-k,k]}$}.			
		\end{array} \right.
	\ees 
\noindent	
Hence, we may say that two biinfinite sequences are close to each other  when their central blocks  agree. We denote the shift map $\phi\colon\mathcal A\sp\Z\to\mathcal A\sp\Z$ as the map which moves the zeroth coordenate one slot to the right, that is $\phi(b)_i=b_{i+1}$. A subshift $X$ is a subspace of $\mathcal A\sp\Z$ which is closed and invariant under $\phi$.
	The  {orbit} of a point $b \in \mathcal A\sp\Z$ is the the set of iterates $\conj{\phi^n(b)}_{n\in\Z}$ and we will denote it as $\orb{b}$.  Given $b\in\mathcal A\sp\Z$ , the subshift generated by $b$ is defined to be $\ovl{\orb{b}}$ and denoted by $X_b$. Notice that when $b$ is pseudoergodic, it follows that $X_b=\mathcal A\sp\Z$; indeed, if $c$ belongs to $\mathcal A\sp\Z$,  then for any $k$ the central block $c_{[-k,k]}$ occurs in $b$, so there is $j_0$ such that  $c_{[-k,k]}=b_{[j_0,j_0+2k]}$. But $b_{[j_0,j_0+2k]}=\phi^{j_0+k}(b)_{[-k,k]}$ so that there is always an element in $\orb{b}$ as close as desired to $c$, proving that $c$ belongs to $X_b$, as wanted.

For $\mathfrak a\in\mathcal A$, we denote its corresponding biinfinite sequence in $\mathcal A\sp\Z$ with boldface 
${\boldsymbol{\mathfrak a}} = (\cdots \mathfrak a\,\mathfrak a\,
\dot{\mathfrak{a}}\,\mathfrak a\cdots)$. The tridiagonal operators 
$A_{\boldsymbol{\mathfrak a}}$ correspondig to the constant sequences $\boldsymbol{\mathfrak a}$ are known as Laurent operators.

\section{Tridiagonal operators for biinfinite sequences of zeroes and ones}

In this section we discuss tridiagonal operators $A_b$ where $b$ is a biinfinite sequence of symbols in the alphabet  $\mathcal A=\{0,1\}$.
Similar arguments can be given for $\mathcal A= \{-1,1\}$ and we will comment about this at the end of this section.
We begin this section by computing the numerical range of the Laurent operators $A\sb\bz$ and $A\sb\bo$; although it is a well known result, we provide a sketch of the proof since we believe it sheds some light in the understanding of the computation of $W(A_b)$, for more general $b$,  which is one of the main result of this section.

\begin{proposition}\label{WAzerosigualD}
	$W\parn{A_\bz}=\D$ and $W\parn{A_\bo}=(-2,2)$.
\end{proposition}
\begin{pf}
For the proof of the inclusions $W\parn{A_\bz}\supset \D$ and $W\parn{A_\bo}\supset (-2,2)$, we will show that for each  $\lm\in\D$ there is $x$ in the unit ball of $\ell^2(\Z)$ such that $\inpr{A\sb\bz x,x}=\lambda$ and $\inpr{A\sb\bo x ,x}=2\mathrm{Re}(\lambda)$. Indeed, we define	
	\bes \label{vprbz}
	x_k= \left\{
	\begin{array}{cc}
		 \sqrt{1-|\lambda|^2}\ \lm^k & \mbox{si}\ k\geq 0\\
		0           & \mbox{si}\ k<0.
	\end{array} \right.
	\ees
It follows that $$\|x\|^2=\sum_{k\in\Z}{|x_k|^2}=(1-|\lambda|^2)\sum_{k\geq 0}{(|\lambda|^2)^k}=1$$

Furthermore 

	\begin{align*}
		\inpr{A_\bz x,x} &=\sum_{k\in\Z}{\ovl{x_k}x_{k+1}}\\
		&=\lm\ (1-|\lambda|^2|)\sum_{k\geq 0}{|\lm|^{2k}}\\
		&=\lm
	\end{align*}
and
	\begin{align*}
		\inpr{A_\bo x,x} &=\sum_{k\in\Z}\left(\ovl{x_k}x_{k+1}+\ovl{x_k}x_{k-1}\right)\\
		&=\inpr{A_\bz x,x}+\ovl{\inpr{A_\bz x,x}}\\
		&=2\mathrm{Re}(\lm),
	\end{align*}
as wanted.

	To prove the other inclusions, first notice that since $A\sb\bo$ is self-adjoint then $W(A\sb\bo)\subset\R$ by Theorem~\ref{teospectrum}. Let $x\in\ell^2(\Z)$ such that  $\norma{x}=1$. Since $x$ and $A\sb\bz x$ are linearly independent,  Cauchy-Schwarz inequality gives $\abs{\inpr{A_\bz x,x}}< \norma{A_\bz x}\cdot\norma{x}=1$ and so $\abs{\inpr{A_\bo x,x}}=\abs{2\mathrm{Re}(\lm)}<2$ so that $W(A\sb\bz)\subset\D$ and  $W(A\sb\bo)\subset (-2,2)$, as was to be proved.
	
\end{pf}

\begin{lemma}\label{4ellipse}
Let $x=(x_j)_{j\in\Z}$ be an element in the unit ball of $\ell^2(\Z)$ and let $b=(b_j)_{j\in\Z}$ be in $\{0,1\}\sp\Z$. Then
\begin{align*}
\inpr{A_b x,x}&= \inpr{A_0 x,x}+\sum_{k\in\Z}{b_k x_{k}\ovl{x_{k+1}}}\nonumber \\
               &= 2\mathrm{Re}(\inpr{A_0x,x})-\sum_{k\in\Z} (1-b_k)x_k\ovl{x_{k+1}}.
\end{align*}
\end{lemma}
\begin{pf}
A direct computation shows
\begin{align}
	\inpr{A_b x,x}&=\sum_{k\in\Z }{(b_{k-1}x_{k-1}+x_{k+1})\ovl{x_k}}\nonumber\\
	&=\inpr{A_\bz  x,x}+\sum_{k\in\Z}{b_k x_{k}\ovl{x_{k+1}}}\label{A_b=}
\end{align}

On the other hand we have
		\begin{align*}
	\inpr{A_\bz  x,x}
              &= \sum_{k\in\Z}{\ovl{x_k}x_{k+1}}\\
	\      &=\sum_{k\in\Z}{(1-b_k)\ovl{x_k}x_{k+1}} + \sum_{k\in\Z}{b_k \ovl{x_k}x_{k+1}}\\
	\end{align*}
and so 
     \begin{align*}
	\sum_{k\in\Z}{b_k {x_k}\ovl{x_{k+1}}} 
&=\ovl{\inpr{A_bz x,x}} - \sum_{k\in\Z}{\left(1-b_k\right){x_k}\ovl{x_{k+1}}}.
	\end{align*}
	
Substituting this last equality in (\ref{A_b=}) we complete the proof.
\end{pf}

\begin{proposition}\label{innerProductInEllipse}
	Let  $x=(x_j)_{j\in\Z}$ be in the unit ball of $\ell^2(\Z)$ an let $b\in\conjz$.  Then the complex number $\inpr{A_bx,x}$ is contained in the interior of the ellipse with major axis equal to 1 and focal points $\inpr{A_\bz x,x}$ and $2\mathrm{Re}(\inpr{A_\bz x,x})$.
\end{proposition}

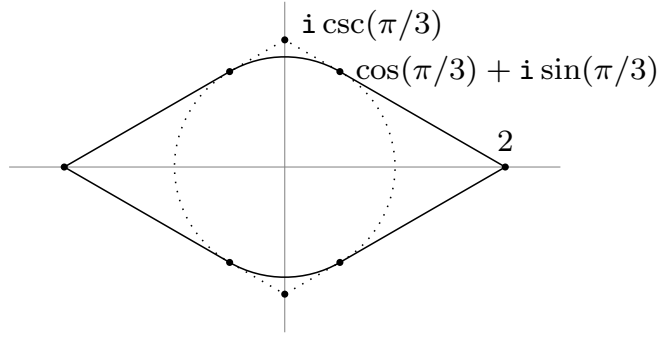
\begin{figure}
\centering
\begin{adjustbox}{min width=0.7 \textwidth}
     \begin{tikzpicture} 
      \draw[help lines] (-2.5,0) -- (2.5,0);
      \draw[help lines] (0,-1.5) -- (0,1.5);
       \draw [black] (0.5,-0.866) -- (2,0) -- (0.5,0.866);
       \draw [black]  (-0.5,0.866) -- (-2,0) -- (-0.5,-0.866);
       \draw [dotted] (0.5,0.866) -- (0,1.154) -- (-0.5,0.866);
       \draw [dotted] (-0.5,-0.866) -- (0,-1.154) -- (0.5,-0.866);
       \draw [dotted] (-0.5,0.866) arc [radius=1, start angle=120, end angle=240];
       \draw [dotted] (0.5,-0.866) arc [radius=1, start angle=300, end angle=420];
	\draw [black] (0.5,0.866) arc [radius=1, start angle=60, end angle=120 ];
	\draw [black] (-0.5,-0.866) arc [radius=1, start angle=240, end angle=300 ]; 
	\draw[fill] (0.5,0.866) circle [radius=0.025];
	\draw[fill] (0.5,-0.866) circle [radius=0.025];
	\draw[fill] (-0.5,0.866) circle [radius=0.025];
	\draw[fill] (-0.5,-0.866) circle [radius=0.025];
      \draw[fill] (2,0) circle [radius=0.025];
      \draw[fill] (-2,0) circle [radius=0.025];
      \draw[fill] (0,1.154) circle [radius=0.025];
      \draw[fill] (0,-1.154) circle [radius=0.025];
	\node [right] at (0.5,0.866) {\tiny $\cos(\pi/3)+\im\sin(\pi/3)$}; 
	\node [above] at (2,0) {\tiny $2$};
	\node [above right] at (0,1) {\tiny $\im\csc(\pi/3)$};  
     \end{tikzpicture}
\end{adjustbox}
\caption{Boundary of  $\Gamma= \conv{W(A\sb\bz)\cup W(A\sb\bo)}$}\label{Gamma}
\end{figure}

\begin{pf}
     For convenience, let us denote by $\Gamma$ the right hand side of the inclusion in the theorem, see Figure~\ref{Gamma}. Since $x$ is in the unit ball of $\ell^2(\Z)$ then the element $y$ in $\ell^2(\Z)$ defined by $y_k=|x_k|$ is also in the unit ball. Therefore, by Proposition~\ref{WAzerosigualD}, we obtain $\inpr{A\sb\bz y,y}<1$. This together with Lemma~\ref{4ellipse} give us
\begin{align*}
	\abs{\inpr{A_bx,x}-\inpr{A_\bz x,x}} & +|\inpr{A_bx,x}-2\mathrm{Re}\inpr{A_\bz x,x}| \\ 
  &=   \abs{\sum_{k\in\Z}b_k {x_k}\ovl{x_{k+1}}}+\abs{-\sum_{k\in\Z}(1-b_k)x_k\ovl{x_{k+1}}}\\
	&\leq \sum_{k\in\Z,\ b_k=1}\abs {\ovl{x_k}x_{k+1}} + \sum_{k\in\Z,\ b_k=0} \abs{{x_k} \ovl{x_{k+1}}}\\
	&=\sum_{k\in\Z}\abs{\ovl{x_k}{x_{k+1}}}\\
       &=\inpr{A\sb\bz y , y }\\
	&<1.
	\end{align*}	
as was to be proved.	
\end{pf}

\begin{definition}
  Let $x=(x_j)_{j\in\Z}$ be in the unit ball of $\ell^2(\Z)$. We will denote by $E_x$ the open convex set limited by the ellipse with focal points at $\inpr{A_0x,x}$ and $2\mathrm{Re}(\inpr{A_0x,x})$ and major axis of length 1.
\end{definition}

\begin{theorem}\label{theoremsubset}
	Let $b\in\conjz$, then 
	\bes
		W(A_b)\subseteq\conv{W(A\sb\bz)\cup W(A\sb\bo)}.
	\ees
\end{theorem}

\begin{pf}
   For convenience, let us denote by $\Gamma$ the right hand side of the inclusion in the theorem, see Figure~\ref{Gamma}.  Using Proposition~\ref{innerProductInEllipse}, it will suffice to show that, for each $x=(x_j)_{j\in\Z}$  in the unit ball of $\ell^2(\Z)$, the corresponding ellipse $E_x$  is contained in $\Gamma$. This is achieved by showing that the bounday of $\Gamma$ lies outside of each $E_x$.  By symmetry, we may only prove the case when $\inpr{A_\bz x,x}$ lies in the first quadrant.  Since the case $\inpr{A_\bz x,x}=0$ implies $E_x$ is the open unit ball, satisfying the desired property, we may assume that $\lambda=\inpr{A_\bz x,x}\not=0$, say $\lambda=r e^{\im \theta}$ with $0<r<1$ and $0\leq \theta\leq \pi/2$. Since the distance from $2\mathrm{Re}(\lambda)$ to the closest boundary line is $\frac{1}{2}\left( 2-2\mathrm{Re}(\lambda) \right)=1-r\cos\theta$,  by {\em Heron's Problem},  the minimum sum of the distances from the boundary line of $\Gamma$ to the points $\lambda$ and $2\mathrm{Re}(\lambda)$ is

\begin{align*}
  \left| \Biggl( 2\mathrm{Re}(\lambda )   +  \left(1  - 
 r\cos\theta\right)\left(1 +\im \sqrt{3} \right)\Biggr) -\lambda \right|\mspace{-200mu}\\
      &=  \abs{ 1+ \im \Biggl(\sqrt{3} - r\left( \sqrt{3}\cos\theta+\sin\theta \right)\Biggr)  } \\
      &= \abs{1 + \im \left(\sqrt{3}-2r\cos(\theta-\pi/6)\right) } \\
      &=\sqrt{4+4(r\cos(\theta-\pi/6))^2-4\sqrt{3}(r\cos(\theta-\pi/6))}
\end{align*}

Furthermore, the above expression, as function of $r\cos(\theta-\pi/6)$, reaches its minimum  value 1 when $r\cos(\theta-\pi/6)=\frac{\sqrt{3}}{2}$. This proves that  the sum of the distances from the boundary line to the points $\lambda$ and $2\mathrm{Re}(\lambda)$ has minimal value 1. Hence it lies outside $E_x$, as desired.

    It remains to prove that the sum of the distances from a point in the curved part of the boundary of $\Gamma$ to $\lambda$ and $2\mathrm{Re}(\lambda)$ is  also at least 1. Por this purpose, let $z=e^{i\varphi}$ with $\pi/3< \varphi\leq \pi/2$ and let $\lambda\sp\prime=\frac{1}{r}\lambda$. We consider two cases, namely,  either $0\leq\theta\leq\varphi$ or $\phi<\theta\leq \pi/2$. In case  $0\leq\theta\leq\varphi$, we look at the triangle with vertices  $z$, $\lambda$,and $\lambda\sp\prime$ and the triangle with vertices $z$, $\lambda$,  and $2\mathrm{Re}(\lambda)$, see Figure~\ref{Fig:Case1}. Since the longest (shortest) side is always opposite the largest (smallest) interior angle, we obtain $|z-\lambda|\geq|\lambda-\lambda\sp\prime|$ and  $|z-2\mathrm{Re}(\lambda)|\geq |\lambda-2\mathrm{Re}(\lambda)|$. Since $|\lambda-2\mathrm{Re}(\lambda)|=|\lambda|$, we obtain 
\begin{align*}
|z-\lambda|+ |z-2 \mathrm{Re}(\lambda)|
&\geq |\lambda-\lambda\sp\prime| + |\lambda| =1.
\end{align*}
as wanted. 

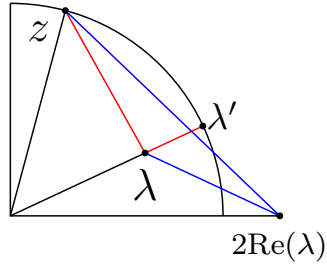
\begin{figure}
\centering
\begin{adjustbox}{min width=0.35 \textwidth}
     \begin{tikzpicture} 
       \draw [black] (0,2) -- (0,0) -- (1.4*1.81,0);
	\draw [black] (2,0) arc [radius=2, start angle=0, end angle= 90];
	\draw (0,0) -- (0.7*1.81,0.7*0.845);
	\draw [red]  (0.5176,1.9318) -- (0.7*1.81,0.7*0.845) -- (1.81,0.845);
	\draw [blue] (0.7*1.81,0.7*0.845) -- (1.4*1.81,0) -- (0.5176,1.9318) ;
	\draw (0,0) -- (0.5176,1.9318);
	\draw[fill] (0.7*1.81,0.7*0.845) circle [radius=0.025];
	\node [below] at (0.7*1.81,0.7*0.845) {$\lambda$};
	\draw[fill] (1.81,0.845) circle [radius=0.025];
	\node [right] at (1.70,0.945) {\footnotesize$\lambda\sp\prime$};
	\draw[fill] (0.5176,1.9318) circle [radius=0.025];
	\node [left] at (0.5176,1.7318) {$z$};
	\draw[fill] (1.4*1.81,0) circle [radius=0.025];
	\node [below] at (1.4*1.81,0) {\tiny$2\mathrm{Re}(\lambda)$};
    \end{tikzpicture}
\end{adjustbox}
\caption{Case $0\leq\theta\leq\varphi$}\label{Fig:Case1}
\end{figure}
	
To complete the proof we now assume $\varphi<\theta\leq \pi/2$.
We use \emph{Ptolemy's inequality} (see Figure~\ref{Fig:Case2}) on the quadrilateral with vertices $\lm$, $z$, $2\mathrm{Re}(z)$ and the origin to obtain 

\begin{align*}
 |z- \lambda|\ |2\mathrm{Re}(\lm)| + |\lm|\  |z-2\mathrm{Re}(\lm)|
       &\geq |\lm - 2\mathrm{Re}(\lm)|
\end{align*}
Using $|2\mathrm{Re}(\lm)|=|2r\cos\theta|\leq |2r\cos(\pi/3)|=r=|\lm|$ and again that $|\lambda-2\mathrm{Re}(\lambda)|=|\lambda|$ , we get

\begin{align*}
|z-\lambda|\ |\lm| + |\lm|\  |z-2\mathrm{Re}(\lambda)|
 &\geq |\lambda - z|\ |2\mathrm{Re}(\lm)| + |\lm|\  |z-2\mathrm{Re}(\lm)|\\
       &\geq |\lambda|
\end{align*}
and so
\begin{align*}
 |z- \lambda| +   |z-2\mathrm{Re}(\lambda)| &\geq 1
\end{align*}
as was to be proved.
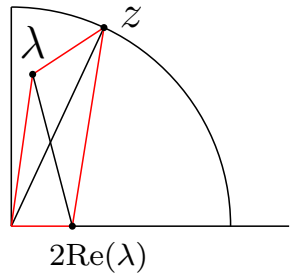
\begin{figure}[ht]
\centering
\begin{adjustbox}{min width=0.3 \textwidth}
\begin{tikzpicture}
\draw [black] (0,0) -- (0.8452,1.8126);
\draw [black] (0.7*0.2783,0.7*1.98) -- (2*0.2783,0);
\draw [black] (0,2) -- (0,0);
\draw [black] (2*0.2783,0) -- (1.4*1.81,0);
\draw [black] (2,0) arc [radius=2, start angle=0, end angle= 90];
\draw [red]  (0,0) -- (2*0.2783,0) -- (0.8452,1.8126) -- (0.7*0.2783,0.7*1.98) 
              -- (0,0);
\draw[fill] (0.7*0.2783,0.7*1.98) circle [radius=0.025];
\node [above] at (0.7*0.2783,0.7*1.98) {$\lambda$};
\draw[fill] (0.8452,1.8126) circle [radius=0.025];
\node [right] at (0.8452,1.9126) {$z$};
\draw[fill]  (2*0.2783,0) circle [radius=0.025];
\node [below] at  (2*0.4,0) {\tiny $2\mathrm{Re}(\lambda)$};
\end{tikzpicture}
\end{adjustbox}
\caption{Case $\varphi<\theta\leq \pi/2$}\label{Fig:Case2}
\end{figure}

\end{pf}

As mentioned in the introduction, the following result might me derived from limit operator techniques when the closure of the numerical range is considered, see e.g. \cite{Hagger}. However, our approach might be of interest since it is elementary and provides the construction of the unitary elements in $\ell^2(\Z)$ correspondig to a given element in the numerical range.

\begin{proposition}\label{InGamma}
Let $b\in\conjz$.
\begin{enumerate}
 \item If the constant biinfinite sequence $\bz$ is in $\ovl{\orb{b}}$ then $$W(A_\bz)=\D\subset W(A_b).$$
\item  If the constant biinfinite sequence $\bo$ is in $\ovl{\orb{b}}$ then $$W(A_\bo)=(-2,2)\subset W(A_b).$$ 
\end{enumerate}
\end{proposition}

\begin{pf}
To prove the proposition, it will suffice to show that  for each $\lm\in\D$ we can provide  elements $x$ and $y$ in the unit ball of $\ell^2(\Z)$ such that $\inpr{A_bx,x}=\lambda$ and $\inpr{A_by,y}=2\mathrm{Re}(\lambda)$.

 	We first observe that the case $\lambda=0$ is taken care of by  choosing $x\in\ell^2(\Z)$ with just one component equal to 1 and all  others equal to zero since then $\inpr{A_bx,x}=0$. So we may assume in what follows that  $\lambda\not =0$.  Let us assume we have $\lm\in\D$, say $\lm=re^{\im\theta}$.  Since  $0<r<1$, then $0<\frac{r+1}{2}<1$. Let  $k_0$ be an integer such that
 	
 	\beq \label{condf}
 		\frac{1-r}{2}>\parn{\frac{r+1}{2}}^{2k_0+1}.
 	\eeq
 	
 	Consider the polynomial function $f(t)=t-t^{2k_0+1}-r$. Note that (\ref{condf})  implies $f(r)=-r^{2k_0+1}<0$ and $f(\frac{r+1}{2})=\frac{1-r}{2}-\parn{\frac{r+1}{2}}^{2k_0+1}>0$. The Intermediate Value Theorem gives us then a $r<t_0<\frac{r+1}{2}<1$ such that
 	
 	\begin{align}
 		t_0-t_0^{2k_0+1}-r=0\nonumber\\
 		1-t_0^{2k_0}=\frac{r}{t_0}. \label{eqtcero}	
 	\end{align}

 	Now, since the constant sequence 0 (resp. 1) is a limit point of the orbit of $b$, there exists an integer  $j_0$ (resp. $l_0$) be such that $b_{[j_0+1,j_0+k_0]}$ (resp. $b_{[l_0+1,l_0+k_0]}$) is a sub-word of zeros (resp. of ones) of $b$ with length $k_0$.   We define  $x=(x_k)_{k\in\Z}$ in $\ell^2(\Z)$  in the following way.  For each $k\in\Z$
 	
 	\bes 
		x_k= \left\{
		\begin{array}{ll}
			\sqrt{1-t_0^2}\ t_0^{k-j_0-1\ }e^{\im\theta(k-j_0)} 
                                 & \text{ if }  j_0<k<j_0+k_0+2\\
			t_0^{k_0+1}				& \mbox{if}\ k=j_0+k_0+3\\
			0            		& \mbox{otherwise.}
		\end{array} \right.
	\ees 
and define  $y_k$ similarly by replacing $j_0$ with $l_0$ in the definition of $x_k$.
 	It follows that
 	
 	\begin{align*}
		\norma{x}^2 
             &=\sum_{k\in\Z}{|x_k|^2}\\ 
		&=t_0^{2(k_0+1)}+\sum_{k=0}^{k_0}(1-t_0^2)t_0^{2k} \\
		&=t_0^{2(k_0+1)}+ (1-t_0^2)\frac{1-t_0^{2(k_0+1)}}{1-t_0^2}\\
             &=1
	\end{align*}
and similarly $\norma{y}=1$.		
	Thus $x$ and $y$  are elements in the unit ball of $\ell^2(\Z)$. Using $b_{j_0+1}=\cdots=b_{j_0+k_0}=0$  we obtain
	
	\begin{align*}
		\inpr{A_b x,x} 
            &= \sum_{k\in\Z}{b_k x_k \ovl{x_{k+1}} + \ovl{x_k}x_{k+1}}\\
		&=\sum_{k=j_0+1}^{j_0+k_0}\ovl{x_k} x_{k+1}\\
	 	&=\sum_{k=j_0+1}^{j_0+k_0} (1-t_0^2)
                    t_0^{2(k - j_0)} t_{0}^{-1}e^{\im\theta}\\
		&= t_0 (1-t_0^{2 k_0}) e^{\im\theta} \\
            &= r e^{\im\theta}\\
            &=\lambda.
	\end{align*}
	
 Similarly 
\begin{align*}
      \inpr{A_by,y}
            &=\ovl{\inpr{A_bx,x}}+\inpr{A_bx,x}\\
            &=\ovl{\lambda}+\lambda\\
            &=2\mathrm{Re}(\lambda)
\end{align*}
as was to be proved.
\end{pf}

 The following theorem is analogous to \cite[Lemma~3.1]{CWChLi10} for $b$ pseudoergodic and the alphabet $\mathcal A=\{-1,1\}$. 
  In case we replace the numerical range with its closure, our result is analogous to \cite[Corollary~14]{Hagger}, where the operator is taken to be pseudoergodic and the union is taken over all periodic operators.

\begin{theorem}\label{theoremequal}
	Let $ b\in\conjz$. If  $\bz,\bo\in X_b$ then
	\bes
		W(A_b)=\conv{W(A\sb\bz)\cup W(A\sb\bo)}.
	\ees
\end{theorem}
\begin{pf}
Since $W(A_b)$ is a convex set, it follows from Proposition~\ref{InGamma} that $W(A\sb\bz)\cup W(A\sb\bo)\subseteq W(A_b)$. The other inclusion is an application of Theorem~\ref{theoremsubset}.
\end{pf}

Our following corollary may also be derived from \cite[Theorem~16]{Hagger}. However, notice we do not requerie the operator $A_b$ to be pseudoergodic and we do not have employed thus far limit operator techniques.

\begin{corollary}
Let $b\in\{0,1\}\sp\Z$. If $\bz,\bo\in X_b$ 
then
\[
\ovl{W(A_b)}=\conv{ \sigma(A_\bz)\cup\sigma(A_\bo)}
\]
\end{corollary}

\begin{pf}
Since $A\sb\bz$ and $A\sb\bo$ are normal operators, Theorem~\ref{teospectrum} gives us that $\conv{\sigma(A\sb\bz)}=\ovl{W(A\sb\bz)}$ and $\conv{\sigma(A\sb\bo)}=\ovl{W(A\sb\bo)}$. Furthermore, since the numerical range of an operator is a convex set, by applying the closure to both sides of the equality in Theorem~\ref{theoremequal}, a straightforward computation gives the desired result.
\end{pf}

Notice that in the proof of Theorem~\ref{InGamma}, by replacing $0$ with $-1$,  it is easy to see that actually $\inpr{A_{\boldsymbol{-1}}x,x}=\ovl{\lambda}$. Hence,  it might be possible to use our techniques to recover some known results for the hopping sign model operators.

\section{Tridiagonal operators associated to subshifts}

For this section, $\mathcal A$  denotes a finite alphabet. 
\begin{proposition}\label{orbitProp}
	Let $b\in\mathcal A\sp\Z$. If $d\in \orb{b}$ then the operators $A_b$ and $A_d$ are unitarily equivalent.
\end{proposition}

\begin{pf}
	Since $d\in \orb{b}$, then exist $k_0\in\Z$ such that $d=\phi^{k} (b)$. We define $S=A_0^{k}$. Note $S$ is an unitary operator. For $j\in\Z$ it follows that
	
	\begin{align*}
		\parn{\parn{S A_b} x}_j&=\parn{A_b x}_{j+k}\\
		&=x_{j+k+1}+b_{j+k-1}x_{j+k-1}\\
		&=\parn{S x}_{j+1}
			+\parn{\phi^{k}( b)}_{j-1}\parn{S x}_{j-1}\\
		&=\parn{Sx}_{j+1}
			+d_{j-1}\parn{S x}_{j-1}\\
		&=\parn{\parn{A_d S }x}_j.
	\end{align*}
		Thus $S^{-1} A_bS= A_d$.
\end{pf}

\begin{lemma}\label{orbitLemma}
If $A$ and $B$ are unitarily equivalent, then they have the same numerical range and spectrum, furthermore, they have the same eigenvalues.
\end{lemma}
\begin{pf}
Since $A$ and $B$ are unitarily equivalent there exists a unitary operator $U$ such that $A=U^{-1}BU$. Let be $\lm\in\C$, it follows that
	\begin{align*} 	
            A-\lm I &= P^{-1}\parn{B-\lm I}P .	
	\end{align*}
	
	Therefore $A-\lm I$ is an 	invertible operator with bounded inverse and dense range if and only if so is $B-\lm I$ is. Hence  $\si\parn{A}=\si\parn{B}$ and $A$ and $B$ have the same eigenvalues.

Finally, since $U^{-1}=U\sp\ast$, it follows that $$\inpr{Ax,x}=\inpr{U^{-1}BUx,x}=\inpr{BUx,Ux}$$
and thus since $U$ is an isometry we conclude $W(A)=W(B)$.
\end{pf}

\begin{theorem} \label{orbit=same}
Let $b\in\mathcal A\sp\Z$, then for every $d\in \orb{b}$  we have
	$	W\parn{A_d}=W\parn{A_b} \mbox{  and }
		\si\parn{A_d}=\si\parn{A_b}$.  
\end{theorem}
\begin{pf}
Follows from Lemma~\ref{orbitLemma} and Proposition~\ref{orbitProp}. 
\end{pf}

The following may be derived from \cite[Lemma~12]{Hagger}; however, we provide a simplified proof here for completeness sake. 

  \begin{lemma}\label{weak}
Let $b,d\in \mathcal A\sp\Z$. If $\phi^{n_k}(b)\to d$ for some sequence of integers $\{n_k\}_{k>0}$ then $A_{\phi^{n_k}(b)}\to A_d$  in the weak operator topology.
 \end{lemma}
\begin{pf}
   Let $m>0$. Since $\phi^{n_k}(b)\to d$, there is $M>0$ such that $d_{[-m,m]}=\phi^{n_k}(b)_{[-m,m]}=b_{[-m+n_k,m+n_k]}$ for all $k>M$.  Since $m$ is arbitrary, this yields $\inpr{A_{\phi^{n_{k}}\left(b\right)}x,y}$ to be as close as wanted to $\inpr{A_dx,y}$ for $k$ sufficiently large, as desired. 
\end{pf}

\begin{proposition}\label{inclusion} 
Let $b\in\mathcal A\sp\Z$. If $d\in X_b$ then
\[
W(A_d)\subset\ovl{W(A_b)}\hspace{1cm}\text{and}\hspace{1cm}
\sigma(A_d)\subset \ovl{W(A_b)}.
\]
\end{proposition}
\begin{pf}
Since $d\in\ovl{\mathrm{Orb}(b)}$, there exists a sequence of integers $\{n_k\}_{k\in\Z}$ such that $\phi^{n_k}(b)\to d$ as $k\to\infty$.  	 Using 	
 Lemma~\ref{weak},  we also have $A_{\phi^{n_k}(b)}\to A_d$  in the weak operator topology. To prove the first inclusion, we let $z\in W(A_d)$. Then there is $x\in\ell^2(\Z)$ with $\|x\|=1$ such that $z=\inpr{A_dx,x}$. Hence 
$\inpr{A_{\phi^{n_k}(b)}x,x}-z=\inpr{(A_{\phi^{n_k}(b)}- A_d)x,x}\to 0$ as $k\to\infty$ and so $\inpr{A_{\phi^{n_k}(b)}x,x}\to z$.  Using Theorem~\ref{orbit=same}, we obtain that $\inpr{A_{\phi^{n_k}(b)}x,x}$ belongs to $W\parn{A_{\phi^{n_k}(b)}}=W\parn{A_b}$, this proves $z\in\ovl{W(A_b)}$, as wanted. The  second inclusion is inmediate from Theorem~\ref{teospectrum}.
\end{pf}

\begin{corollary}
    Let $b\in\mathcal A\sp\Z$. Then 
\[
\bigcup_{d\in X_b}\ovl{W(A_d)}=\ovl{W(A_b)}.
\]
\end{corollary}
\begin{pf}
    The nontrivial inclusion is an application of Proposition~\ref{inclusion}
\end{pf}

The following corollary is analogous to \cite[Theorem~16]{Hagger}, in fact, we rely in its proof.

\begin{corollary}\label{igualdad}
Let $b\in\mathcal A\sp\Z$. If for all $\mathfrak a\in\mathcal A$ the constant sequence $\boldsymbol{\mathfrak{a}}$ belongs to  $X_b$ then
\[
   \ovl{W(A_b)}=\conv{\bigcup_{\mathfrak a\in\mathcal A} \ovl{W(A_{\boldsymbol{\mathfrak{a}}})}}=\conv{\bigcup_{\mathfrak{a}\in\mathcal A}\sigma(A_{\boldsymbol{\mathfrak{a}}})}.
\]
\end{corollary}
\begin{pf}
     The equality in the right hand side follows from Theorem~\ref{teospectrum} as the Laurent operator $A_{\boldsymbol{\mathfrak{a}}}$ is self-adjoint. For the left hand side equality, the inclusion ``$\supset$'' follows from an application of Theorem~\ref{orbit=same} by taking closures and using that the closure of a convex set is convex. The inclusion ``$\subset$'' follows from the proof  of \cite[Theorem~16]{Hagger} by taking $U_{-1}=\mathcal A$, $U_0=\{0\}$ and $U_1=\{1\}$ and noticing that the hypothesis on $A=A_b$ to be pseudoergodic is not required in that part of the proof.
\end{pf}

The following corollary is somewhat more general than a particular case of \cite[Corollary~17]{Hagger}, as applies for not necessary pseudoergodic sequences.

\begin{corollary}\label{unionOfAll}
 Let $b\in\mathcal A\sp\Z$. If ${\boldsymbol{\mathfrak a}} \in X_b$ for every $\mathfrak a\in\mathcal A$, then 
\[
\ovl{W(A_b)}=\conv{\sigma(A_b)}.
\]
\end{corollary}
\begin{pf}
Since $b\in X_b$, using Proposition~\ref{inclusion} and Theorem~\ref{teospectrum} we conclude the ``$\supset$'' inclusion. For the other inclusion we observe that  the hypothesis ensure that Laurent operators $A_{\boldsymbol{\mathfrak a}}$ are limit operators of $A_b$ (see e.g. \cite[Section~2]{CWChLi13} for an introduction of this notion)  and furthermore that $A_b$ is actually a limit operator of itself.  We  may apply now \cite[Theorem~2.1]{CWChLi13} to obtain   
\[
     \bigcup_{\mathfrak a\in\mathcal A} {\sigma(A_{\boldsymbol{\mathfrak{a}}})}
\subset \sigma(A_b) 
\]
This together with Corollary~\ref{igualdad} complete the proof.
\end{pf}

Since subshifts may not contain periodic points at all, e.g.~the Sturmian subhifts with irrational frequency, it is false in general to restrict the union on left hand side of Corollary~\ref{unionOfAll} by periodic points. 
However, if $X_b$ is of finite type, the periodic points are dense \cite[Exercise~6.1.12]{LinMar95}; such is the case when $b$ is pseudoergodic. We conjecture that when $X_b$ is of finite type the union of the spectrums of the tridiagonal operators corresponding to periodic sequences in $X_b$ is dense in the spectrum of $A_b$. For the case of pseudoergodic sequences this was conjectured in \cite{CWChLi13} and proved recently in \cite{HaggerJFA}.


\end{document}